\documentclass[11pt,draft]{amsart}
\usepackage{setspace,units}
\usepackage{bbm,latexsym,mathrsfs,amscd}

\newtheorem{theorem}{Theorem}[section]
\newtheorem{lemma}[theorem]{Lemma}

\theoremstyle{definition}
\newtheorem{definition}[theorem]{Definition}
\newtheorem{example}[theorem]{Example}

\numberwithin{equation}{section}

\RequirePackage{amsmath}
\RequirePackage{amssymb}
\RequirePackage{amsthm}
 \RequirePackage{epsfig}

\begin{document}

\date{}
\title[Arithmetical operations on infinite decimals]
{Defining arithmetical operations on infinite decimals}

\author{Nicolas Fardin and Liangpan Li}

\address{Lgt De Lattre de Tassigny, 85000 La Roche-sur-Yon, France}
\email{Nicolas.Fardin@ac-nantes.fr}

\address{Department of Mathematical Sciences,
Loughborough University, LE11 3TU, UK}
 \email{L.Li@lboro.ac.uk, liliangpan@gmail.com}

\subjclass[2010]{97F50}

\keywords{real number system, decimal representation}


\begin{abstract}
Completing Loo-Keng Hua's  approach to the real number system pioneered in 1962,
this paper  defines arithmetical operations directly on infinite decimals without appealing to any ordering structure.
 Therefore, the widespread belief  that  there exists an algorithm
for determining the digits of the product of two real numbers  in terms of finite pieces of their decimal strings
 is  essentially confirmed.
\end{abstract}


\maketitle

The real number system (RNS) was first constructed via partitions of $\mathbb{Q}$ by Dedekind in 1872.
In the same year, Cantor provided   a second approach
in terms of Cauchy sequences of rational numbers. Since then, plenty of attempts have been made to construct the RNS from
various perspectives (see e.g. \cite{Weiss}).
In a series of  articles  \cite{FardinLi,Gowers,Klazar,Li}, the soon-derived least upper bound property
of several slightly different ambient spaces such as $\mathbb{Z}\times\mathbb{Z}_{10}^{\mathbb{N}}$ has been used to define
arithmetical operations. Here $\mathbb{Z}_{10}$ denotes the set $\{0,1,2,\ldots,9\}$.
Completing Hua's  approach to the RNS pioneered in 1962 (\cite{Hua}),
this paper  defines these operations directly on infinite decimals. Therefore,
the widespread belief (see e.g. \cite{Katz}) that  there exists an algorithm
for determining the digits of the product of two real numbers  in terms of finite pieces of their decimal strings
 is essentially confirmed.

Our ambient space is (see e.g. \cite{Courant,Kodaira})
\[\mathbb{R}=\{a_0.a_1a_2a_3\cdots\in\mathbb{Z}\times\mathbb{Z}_{10}^{\mathbb{N}}: a_k<9\ \mbox{for infinitely many}\ k\}.\]
As usual, an element $x=a_0.a_1a_2a_3\cdots$  is said to be terminating  if there exists a non-negative integer $m$ such that $a_k=0$ for $k>m$.
In this case,  write $x=a_0.a_1a_2\cdots a_m$ for simplicity.
Defining addition and multiplication on the collection of all terminating decimals is rather standard.
For example,
\[1.2\times(-3).4+5.6=-(1.2\times2.6)+5.6=-(3.12)+5.6=(-4).88+5.6=2.48.\]

Hua's idea of defining addition on $\mathbb{R}$ is as follows (see also \cite{Richman}).
For any element $x=a_0.a_1a_2a_3\cdots$ and any non-negative integer $k$, denote $\theta_k(x)=a_k$, the $k$-th digit of $x$, and
$x_k=a_0.a_1a_2\cdots a_k$,  the rational truncation of $x$ up to the $k$-th digit.
Given $x,y\in\mathbb{R}$, write $x=x_k+\epsilon_k$, $y=y_k+\delta_k$, where the tails
$\epsilon_k$ and $\delta_k$ lie in $[0,10^{-k})$.
Note for $k\geq1$,
\[x+y=x_k+y_k+\epsilon_k+\delta_k=(x_k+y_k)_{k-1}+{\theta_k(x_k+y_k)}\cdot{10^{-k}}+\epsilon_k+\delta_k.\]
So if $\theta_k(x_k+y_k)\leq8$ for some $k\geq1$, then
\[0\leq (x+y)-(x_k+y_k)_{k-1}<10^{-(k-1)},\]
which implies $(x+y)_{k-1}=(x_k+y_k)_{k-1}$.
 Consequently, if there are infinitely many positive integers $k$ such that $\theta_k(x_k+y_k)\leq8$, then $x+y$ is determined iteratively  by this procedure.
One needs also to analyze the case of $\theta_k(x_k+y_k)=9$ for large enough $k$, but fortunately there
 is a natural definition of $x+y$ under this condition.

  In much the same way, we propose a similar definition of multiplication that was not studied in \cite{Hua}.
To verify that these operations form a field, we follow the arguments in \cite{FardinLi}.

\section{Definitions and examples}

We mainly collect and  propose definitions of arithmetical operations in this section, and
will justify them in the next one.

\begin{definition}[addition, \cite{Hua}]\label{addition} Let $x=a_0.a_1a_2a_3\cdots$, $y=b_0.b_1b_2b_3\cdots$ be elements of $\mathbb{R}$.\\
Case 1: Suppose there exists a non-negative integer $m$ such that $a_k+b_k=9$ for $k>m$. Then define\footnote{Here $10^{-m}$ is identified with the terminating decimal $0.\underbrace{00\cdots01}_{\mbox{$m$ digits}}$.}
\[x+y=x_m+y_m+10^{-m}.\]
Case 2: Suppose there exists a sequence of positive integers $k_1<k_2<k_3<\cdots$ such that
$a_{k_i}+b_{k_i}\neq9$ for $i\in\mathbb{N}$. Then $x+y$ is defined by setting
\[(x+y)_{k_i-1}=(x_{k_i}+y_{k_i})_{k_i-1}\ \ \ (i\in\mathbb{N}).\]
\end{definition}

\begin{definition}[additive inverse,  \cite{FardinLi,Gowers,Kodaira,Li}]\label{addition inverse} Let $x=a_0.a_1a_2a_3\cdots$ be an element of $\mathbb{R}$.\\
Case 1: Suppose there exists a non-negative integer $m$ such that $a_k=0$ for $k>m$. Then define
$-x=-x_m.$\\
Case 2: Suppose there exists a sequence of positive integers $k_1<k_2<k_3<\cdots$ such that
$a_{k_i}>0$ for $i\in\mathbb{N}$. Then define
\[-x=(-1-a_0).(9-a_1)(9-a_2)(9-a_3)\cdots.\]
\end{definition}

\begin{definition}[multiplication]\label{multiplication}  Let $x,y$ be   elements of $\mathbb{R}$.\\
 (1) Suppose $x,y$ are non-negative\footnote{An element $x=a_0.a_1a_2a_3\cdots$ is non-negative or negative if $a_0\geq0$ or $a_0<0$. $x$ is  positive if
 it is non-negative and not zero. Similarly, $x>1$ if $a_0\geq1$ and $x\neq1$, $x\leq1$ if $x>1$ does not hold.}. Fix a non-negative integer $s$ such that $x+y\leq10^{s}$.\\
Case 1: Suppose there exists a non-negative integer $m$  such that $\theta_{k}(x_{k+s}y_{k+s})=9$ for $k>m$. Then define
\[xy=(x_{m+s}y_{m+s})_m+10^{-m}.\]
Case 2: Suppose there exists a sequence of positive integers $k_1<k_2<k_3<\cdots$ such that
$\theta_{k_i}(x_{k_i+s}y_{k_i+s})\neq9$ for $i\in\mathbb{N}$. Then $xy$ is defined by setting
\[(xy)_{k_i-1}=(x_{k_i+s}y_{k_i+s})_{k_i-1} \ \ \ (i\in\mathbb{N}).\]
(2) Suppose $x,y$ are negative. Then define
$xy=(-x)(-y).$\\
(3) Suppose only one of $x$ and $y$ is negative. Then define $xy=-(x(-y))$.
\end{definition}

\begin{definition}[reciprocal]\label{quotient} Let $x$ be a non-zero element of $\mathbb{R}$.\\ (1) Suppose $x$ is positive.
Choose the unique element  $y$ of $\mathbb{R}$ that satisfies
\[xy_k\leq1<x(y_k+10^{-k})\ \ \ (k=0,1,2,\ldots).\]
Then define $x^{-1}=y$.\\ (2) Suppose $x$ is negative. Then define $x^{-1}=-(-x)^{-1}$.
\end{definition}

\begin{example} Let $x=a_0.a_1a_2a_3\cdots=0.777777\cdots$, $y=b_0.b_1b_2b_3\cdots=0.232323\cdots$, and denote $x+y=c_0.c_1c_2c_3\cdots$. Note that
$a_2+b_2=10\neq9$ and $x_2+y_2=1.00$. According to Definition \ref{addition}, $x+y=1.0c_2c_3c_4\cdots$. Similarly, from
$a_4+b_4=10\neq9$ and $x_4+y_4=1.0100$ one can deduce
$x+y=1.010c_4c_5c_6\cdots$. Continuing in this way yields $x+y=1.010101\cdots$.
\end{example}

\begin{example}[$\sqrt{2}$]\label{example26}
Choose the unique positive element  $x$ of $\mathbb{Z}\times\mathbb{Z}_{10}^{\mathbb{N}}$ that satisfies
\[x_k^2< 2<(x_k+10^{-k})^2\ \ \ (k=0,1,2,\ldots),\]
that is, setting $1^2<2<2^2$, $1.4^2<2<1.5^2$, $1.41^2<2<1.42^2$, $1.414^2<2<1.415^2$, $\cdots$, yields
\[x=1.41421356237309504880168872420969807856967187537694\cdots.\]
First, we claim $x\in\mathbb{R}$.
If this is not true, then $x$ is of the form $x=a_0.a_1a_2\cdots a_m999\cdots$. So for  $k\geq m$,
\[x_k^2<2<(x_k+10^{-k})^2=(x_m+10^{-m})^2,\]
which yields
\[0<(x_m+10^{-m})^2-2\leq(x_k+10^{-k})^2-x_k^2\leq 5\cdot 10^{-k}.\]
This is absurd as the  positive element $(x_m+10^{-m})^2-2$ could not be bounded from above by $5\cdot10^{-k}$
when $k$ is large enough. Next, we show that $x^2=2$. 
 Obviously, $x+x\leq2+2=4$, so we take $s=1$ in Definition \ref{multiplication}. Note that
\begin{align*}
2>x_{k+1}^2&=2-(2-x_{k+1}^2)\\
&\geq 2-((x_{k+1}+10^{-k-1})^2-x_{k+1}^2)\\
&\geq2-5\cdot10^{-k-1},
\end{align*}
which yields $(x_{k+1}^2)_{k}=2-10^{-k}=1.\underbrace{999\cdots9}_{\mbox{$k$ digits}},$ and consequently, $\theta_k(x_{k+1}^2)=9$ for all $k>0$.
So taking $m=0$ in Definition \ref{multiplication}, we get
\[x^2=(x_1x_1)_0+1=(1.4\times1.4)_0+1=(1.96)_0+1=1+1=2.\]
\end{example}

\section{Justification of definitions}

\subsection{Justification of addition.} We follow the notations and assumptions in Definition \ref{addition}.
In Case 1, one can easily check that
\begin{equation}\label{F41}x_k+y_k+10^{-k}=x_m+y_m+10^{-m}\end{equation}
for $k>m$. So the definition is independent of the choices of $m$.
In Case 2, we first claim
\begin{equation}\label{F444}(x_{n}+y_{n})_{k_i-1}=(x_{k_i}+y_{k_i})_{k_i-1}\end{equation}
for $n>k_i$. To see this, note
\begin{align*}
x_n+y_n&=x_{k_i}+y_{k_i}+(x_n-x_{k_i})+(y_n-y_{k_i})\\
&=(x_{k_i}+y_{k_i})_{k_i-1}+{\theta_{k_i}(x_{k_i}+y_{k_i})}\cdot{10^{-k_i}}+(x_n-x_{k_i})+(y_n-y_{k_i}).\end{align*}
Since the assumption $a_{k_i}+b_{k_i}\neq9$ is equivalent to $\theta_{k_i}(x_{k_i}+y_{k_i})\leq8$, we get
\[0\leq (x_n+y_n)-(x_{k_i}+y_{k_i})_{k_i-1}<10^{-(k_i-1)},\]
which proves the claim. Consequently, $x+y$ is defined as an element of $\mathbb{Z}\times\mathbb{Z}_{10}^{\mathbb{N}}$.
If $x+y$ is not an element of $\mathbb{R}$, say for example $x+y=c_0.c_1c_2c_3\cdots$ with $c_k=9$ for $k$ bigger than or equal to some $s\in\mathbb{N}$, we then assume without loss of generality that $\theta_s(x_s+y_s)\leq8$. Fixing an $l>s$
so that $x_n\leq x_s+10^{-s}-10^{-l}$ and $y_n\leq y_s+10^{-s}-10^{-l}$ for $n\geq s$, we get
\[x_n+y_n\leq(x_s+y_s+2\cdot10^{-s})-2\cdot10^{-l}\leq c_0.c_1c_2\cdots c_{s-1}99\cdots9\ \ \ (n\geq s),\]
where the last digit $9$ is in the $l$-th decimal place. This is absurd if we choose a large enough $n$ with $\theta_n(x_n+y_n)\leq8$.
Therefore, $x+y$ is an element of $\mathbb{R}$.


\subsection{Justification of multiplication.}
Let $x,y$ be non-negative elements of $\mathbb{R}$. Fix a non-negative integer $s$ such that $x+y\leq10^{s}$\\
Case 1: Suppose there exists a non-negative integer $m$  such that $\theta_{k}(x_{k+s}y_{k+s})=9$ for $k>m$.
First, we claim that
 \begin{equation}\label{F44}(x_{n+s}y_{n+s})_m=(x_{m+s}y_{m+s})_m\end{equation}
 for $n>m$. To verify (\ref{F44}), it suffices to do so for $n=m+1$,  and suppose this is the case.
Then
\begin{align*}
x_{n+s}y_{n+s}&=x_{m+s}y_{m+s}+(x_{n+s}-x_{m+s})y_{m+s}+x_{n+s}(y_{n+s}-y_{m+s})\\
&\leq x_{m+s}y_{m+s}+(x_{n+s}+y_{n+s})\cdot\frac{9}{10^{n+s}}\\
&\leq x_{m+s}y_{m+s}+\frac{9}{10^{n}}.\end{align*}
 Considering $\theta_n(x_{n+s}y_{n+s})=9$ and $\theta_n(x_{m+s}y_{m+s})\leq9$, we get
\begin{align*}
(x_{n+s}y_{n+s})_m+\frac{9}{10^n}&=(x_{n+s}y_{n+s})_n\leq (x_{m+s}y_{m+s}+\frac{9}{10^{n}})_n\\
&\leq(x_{m+s}y_{m+s})_m+\frac{9}{10^n}+\frac{9}{10^n}
\end{align*}
 which implies
 \[0\leq(x_{n+s}y_{n+s})_m-(x_{m+s}y_{m+s})_m\leq\frac{9}{10^n}<\frac{1}{10^m}.\]
 This proves the claim (\ref{F44}).
 Next, we claim that
\begin{equation}\label{F43}(x_{n+s}y_{n+s})_n+10^{-n}=(x_{m+s}y_{m+s})_m+10^{-m}\end{equation}
 for $n>m$.  To verify (\ref{F43}), it suffices to do so for $n=m+1$,  and suppose this is the case.
  Recall $\theta_n(x_{n+s}y_{n+s})=9$, so (\ref{F43}) is equivalent to (\ref{F44}). Therefore,
 the  definition is independent of the choices of $m$. On the other hand, it follows from (\ref{F44}) that
  \begin{equation}(x_{n+s}y_{n+s})_m+10^{-m}=(x_{m+s}y_{m+s})_m+10^{-m}.\end{equation}
 So  the  definition is also independent of the choices of $s$.\\
 Case 2: Suppose there exists a sequence of positive integers $k_1<k_2<k_3<\cdots$ such that
$\theta_{k_i}(x_{k_i+s}y_{k_i+s})\neq9$ for $i\in\mathbb{N}$.
 We claim that
\begin{equation}\label{F555}(x_ny_n)_{k_i-1}=(x_{k_i+s}y_{k_i+s})_{k_i-1}\end{equation}
for $n>k_i+s$. Similar to the verification of the previous case, one gets
\begin{equation}\label{F45}x_ny_n=x_{k_i+s}y_{k_i+s}+\gamma_n\end{equation}
with $0\leq\gamma_n<\frac{1}{10^{k_i}}$.
Considering $\theta_{k_i}(x_{k_i+s}y_{k_i+s})\leq8$, we can write
\begin{equation}\label{F46}x_{k_i+s}y_{k_i+s}=(x_{k_i+s}y_{k_i+s})_{k_i-1}+\epsilon_n\end{equation}
with $0\leq\epsilon_n\leq\frac{9}{10^{k_i}}$. Combining (\ref{F45}) and (\ref{F46}) yields
\[0\leq x_ny_n-(x_{k_i+s}y_{k_i+s})_{k_i-1}<\frac{1}{10^{k_i-1}},\]
which proves the claim (\ref{F555}).
Consequently, $xy$ is defined as an element of $\mathbb{Z}\times\mathbb{Z}_{10}^{\mathbb{N}}$.
To finish the justification, one needs to show that
 $xy\in\mathbb{R}$, which is left as an exercise for
interested readers.

\subsection{Inverses}
We leave the justification of Definition \ref{addition inverse} to interested readers.
Given any $x\in\mathbb{R}$, it is easy to check that $x+0=0+x=x$, $x\times 1=1\times x=x$, and $x+(-x)=0$.
Here $x\times y$ means as usual the product between $x$ and $y$.
So 0 is the additive unit, 1 is the multiplicative unit, and $-x$ is the additive inverse of $x$.
Similar to the proof of $\sqrt{2}\in\mathbb{R}$ in Example \ref{example26},
one can show that the reciprocal of a positive element defined by Definition \ref{quotient}
belongs to $\mathbb{R}$.
 Given a positive element $x=a_0.a_1a_2a_3\cdots\in\mathbb{R}$, let $y=b_0.b_1b_2b_3\cdots$ be the reciprocal of $x$ defined by Definition \ref{quotient}.
 Obviously, $y$ is positive. In the following we explain how to derive $xy=1$.
If $x=1$,  then $y=1$ and thus we have nothing to do. So we can assume $x\neq1$, which implies that $x_ky_k$ is strictly less than 1 for $k\geq0$. Then
\begin{align*}
1>x_ky_k=1-(1-x_ky_k)\geq1-(x(y_k+10^{-k})-x_ky_k)\geq1-\frac{a_0+b_0+2}{10^k},
\end{align*}
from which one can deduce $xy=1$.


\section{Arithmetical laws}

To establish various arithmetical laws, we prepare the following three lemmas.

\begin{lemma}\label{lemma1}
If $x\neq y$, then  there exists an $l\in\mathbb{N}$ such that  $|x_k- y_k|\geq 10^{-l}$ for  $k> l$.
\end{lemma}

\begin{lemma}\label{lemma2}
Let $x,y$ be elements of $\mathbb{R}$. Then $|(x+y)_k-x_k-y_k|\leq 4\cdot 10^{-k}$ for all $k$.
\end{lemma}

\begin{lemma}\label{lemma3}
Let $x,y$ be non-negative elements of $\mathbb{R}$. Then $|(xy)_k-x_ky_k|\leq M\cdot 10^{-k}$ for all $k$, where $M$ is a positive integer depending only on $x$ and $y$.
\end{lemma}

A proof of Lemma \ref{lemma1} is as follows. Assume without loss of generality that
\[x=a_0.a_1a_2a_3\cdots<y=b_0.b_1b_2b_3\cdots.\] Take first a non-negative integer $m$ such that
$x_m<y_m$, then a positive integer $l>m$ so that $a_l\leq8$. For $k>l$, we have
\begin{align*}
y_k-x_k&\geq y_m-\Big(x_m+\Big(\sum_{i=m+1}^l\frac{a_i}{10^i}\Big)+10^{-l}\Big)
\geq y_m-\Big(x_m+\sum_{i=m+1}^l\frac{9}{10^i}\Big)\\ &=(y_m-x_m-10^{-m})+10^{-l}\geq 10^{-l},
\end{align*}
which finishes the proof.
To prove Lemma \ref{lemma2}, we trace the justification of Definition \ref{addition}, and thus
 follow the notations and assumptions therein. In Case 1, one has
\[x+y=x_k+y_k+10^{-k}\ \ \ (k>m),\]
which implies
\[|(x+y)_k-x_k-y_k|=|(x+y)-x_k-y_k|=10^{-k}\]
for $k>m$. In  Case 2, $|(x+y)_{k_i-1}-x_{k_i-1}-y_{k_i-1}|$ is bounded from above by
\[
|(x+y)_{k_i-1}-x_{k_i-1}-y_{k_i-1}|\leq|(x_{k_i}+y_{k_i})_{k_i-1}-(x_{k_i}+y_{k_i})|+|x_{k_i}-x_{k_i-1}|+|y_{k_i}-y_{k_i-1}|,\]
which is less than
$3\cdot10^{-(k_i-1)}.$
Therefore, no matter which case happens, there exist infinite many integers $k$ so that
\begin{equation}\label{F42}|(x+y)_k-x_k-y_k|\leq 3\cdot10^{-k}.\end{equation}
Let $q$ be an arbitrary integer, and let $k>q$ be such that (\ref{F42}) holds.
Then
\[|(x+y)_q-x_q-y_q|\leq|(x+y)_k-x_k-y_k|+3\cdot10^{-q}\leq4\cdot10^{-q},\]
which proves Lemma \ref{lemma2}.  Lemma \ref{lemma3} can be derived similarly, so its proof is omitted.

\textbf{Commutative laws}: $x+y=y+x$, $xy=yx$.

These laws are self-evident.

\textbf{Associative laws}: $(x+y)+z=x+(y+z)$, $(xy)z=x(yz)$.

It follows from Lemma \ref{lemma2} that
\[|((x+y)+z)_k-x_k-y_k-z_k|=|((x+y)+z)_k-(x+y)_k-z_k+(x+y)_k-x_k-y_k|\leq 8\cdot 10^{-k}.\]
Similarly,
\[|(x+(y+z))_k-x_k-y_k-z_k|\leq 8\cdot 10^{-k}.\]
So
\[|((x+y)+z)_k-(x+(y+z))_k|\leq 16\cdot 10^{-k}.\]
If $(x+y)+z$ and $x+(y+z)$ are not the same, then there exists an  $l\in\mathbb{N}$ such that
\[|((x+y)+z)_k-(x+(y+z))_k|\geq10^{-l}\]
for $k>l$. Consequently, $10^{-l}\leq 16\cdot 10^{-k}$ for $k>l$, which is absurd if we let $k=l+2$. This proves the associative law for addition.
In much the same way, one can establish the associative law for multiplication between three non-negative elements.
The general case is left as an exercise for interested readers.

\textbf{Distributive law}: $x(y+z)=xy+xz$.\\
Case 1: Suppose $x, y,z$ are non-negative. One can provide a proof that is  similar to that of the associative law for addition.\\
Case 2: Suppose $y$ and $z$ are of the same sign. Then the law follows from Case 1.\\
Case 3: Suppose  $y$ and $z$ are not of the same sign. We can assume  without loss of generality that $y+z$, $-y$, and $z$ are of the same sign.  According to Case 2,
$x(y+z)+x(-y)=xz,$
which yields $x(y+z)=xy+xz$.

To conclude, $(\mathbb{R},+,\times)$ is a field.

\end{document}